\documentclass[11pt]{article}
\usepackage{amssymb,latexsym,amsfonts,amsthm,amsmath}
\usepackage{enumerate, epsfig}

\newcommand{\tC}{\tilde{C}}
\newcommand{\tE}{\tilde{E}}

\newcommand{\tgi}{\tilde{g}^{-1}}
\newcommand{\D}{\mathcal{D}}
\newcommand{\PP}{\mathcal{P}}
\newcommand{\U}{\mathcal{U}}
\newcommand{\0}{\Omega_{0,1}}
\numberwithin{equation}{section}

\theoremstyle{plain}
\newtheorem{theorem}{Theorem}
\newtheorem{corollary}[theorem]{Corollary}
\newtheorem{lemma}[theorem]{Lemma}
\newtheorem{proposition}[theorem]{Proposition}

\theoremstyle{definition}
\newtheorem{definition}[theorem]{Definition}

\newtheorem{example}[theorem]{Example}

\theoremstyle{remark}

\begin{document}
\title{Zeros of the Derivatives of Faber Polynomials Associated with a Universal Covering Map}
\author{Byung-Geun Oh\thanks{Supported by NSF Grant DMS-0244421}\\
        Department of Mathematics\\
        Purdue University\\
        West Lafayette, IN 47907}
\date{February 18, 2004}
\maketitle

\begin{abstract}
For a compact set $E \subset \mathbb{C}$ containing more than two points, 
we study asymptotic behavior of normalized zero counting measures $\{ \mu_k \}$
of the derivatives of Faber polynomials associated with $E$.
For example if $E$ has empty interior, we prove that $\{ \mu_k \}$ converges in the weak-star topology
to a measure whose support is the boundary of $g(\mathcal{D})$, where 
$g : \{ |z| > r \}\cup \{\infty\} \to \overline{\mathbb{C}} \backslash E$ is a universal covering
map such that $g(\infty) = \infty$ and $\mathcal{D}$ is the Dirichlet domain associated 
with $g$ and centered at $\infty$. 

Our results are counterparts of those of Kuijlaars and Saff (1995)
on zeros of Faber polynomials.
\end{abstract}

\section{Introduction}\label{intro}
 Let $g$ be a function which is holomorphic in a neighborhood of infinity
such that its Laurent series is of the form
\begin{equation}\label{laurent}
 g(w) = w + \sum_{k=0}^\infty b_k w^{-k}.
\end{equation}
\emph{Faber polynomials} $F_k$ associated with $g$ are defined by the
generating function
\begin{equation}\label{genf}
  \frac{g'(w)}{g(w)-z} = \sum_{k=0}^\infty \frac{F_k (z)}{w^{k+1}},
\end{equation}
and the \emph{normalized} derivatives of Faber polynomials, 
$P_k :=  (F_{k+1})'/(k+1)$ for $k=0,1,2,\ldots$, satisfy the equation
\begin{equation}\label{genp}
  \frac{1}{g(w)-z} = \sum_{k=0}^\infty \frac{P_k (z)}{w^{k+1}}.
\end{equation}
For every $k$, both $F_k$ and $P_k$ are monic polynomials of degree $k$.
To study more about Faber polynomials and their derivatives, see for example
\cite{HE}, \cite{AES} and \cite{U3}.

Fore each $k$, let $\nu_k$  be the normalized zero counting measure of $F_k$; i.e., 
\[
\nu_k = (2 \pi k)^{-1} \Delta (\log |F_k|),
\] 
where $\Delta$ represents for
the generalized Laplacian. Similarly, we denote by
$\mu_k :=  (2 \pi k)^{-1} \Delta (\log |P_k|)$ the normalized zero counting measure of $P_k$. Kuijlaars and Saff studied
in \cite{KS} the limit behavior of the measures $\{ \nu_k \}$, and we concern  
the limit behavior of $\{ \mu_k \}$. Especially we are interested in the
case when $g$ is a universal covering map. 

Suppose that $E \subset \mathbb{C}$ is a compact set containing more than two points such that
$\Omega := \overline{\mathbb{C}} \backslash E$ is connected (but not necessarily simply-connected).
By the Uniformization theorem (cf. \cite{Ahl}, Chap.~10), 
there exists a unique number $r = r(E) > 0$ and a unique
normalized universal covering map $g: \Lambda (r) := \overline{\mathbb{C}} \backslash \{|w| \leq r \} \to \Omega$
which has a Laurent expansion of the form \eqref{laurent} at infinity.
In this case Faber polynomials $\{ F_k \}$, 
and their normalized derivatives $\{ P_k \}$, associated with $g$ are also called Faber polynomials,
or the normalized derivatives of Faber polynomials, respectively, associated with $E$. 

Since $E$ contains more than two points, the domain $\Omega = \overline{\mathbb{C}} \backslash E$ carries the unique 
hyperbolic (Poincar\'{e}) metric with constant curvature ($\equiv -1$),
and we define $\tE$ as the union of $\partial E$ and
the points in $\Omega$ which have more than one shortest curve
to $\infty$ with respect to this metric. 

\begin{example}
If $\overline{\mathbb{C}} \backslash E$ is simply-connected, that is, $E$ is connected,
then  $\tE = \partial E$.
\end{example} 

\begin{example}\label{example}
Suppose $E$ is a compact set consisting of three points. Then
$\tE$ is either a topological tripod or a line segment joining
points in $E$. See Proposition~\ref{prop}. For example, if
$E = \{1, \eta, \eta^2 \}$ where $\eta = \exp (2 \pi i /3)$ is a third root of unity,
then $\tE = \{ \eta^j t : j = 0,1,2 \mbox{ and } 0 \leq t \leq 1 \}$ (Figure~1).
If $E = \{ -1,0,1 \}$, then $\tE = \{ z \in \mathbb{R} : -1 \leq z \leq 1 \}$ (Figure~2).
\end{example}

\bigskip

\begin{center}
\begin{picture}(0,0)%
\includegraphics{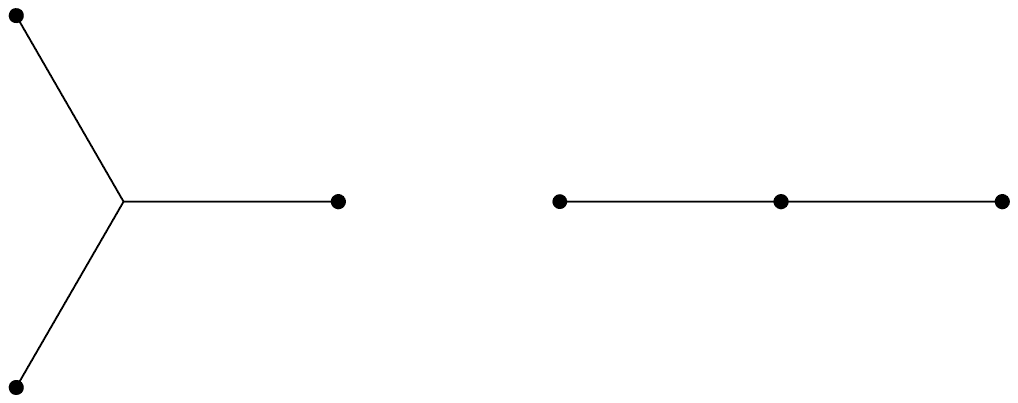}%
\end{picture}%
\setlength{\unitlength}{3947sp}%
\begingroup\makeatletter\ifx\SetFigFont\undefined%
\gdef\SetFigFont#1#2#3#4#5{%
  \reset@font\fontsize{#1}{#2pt}%
  \fontfamily{#3}\fontseries{#4}\fontshape{#5}%
  \selectfont}%
\fi\endgroup%
\begin{picture}(4970,2473)(376,-2819)
\put(676,-2311){\makebox(0,0)[lb]{\smash{\SetFigFont{12}{14.4}{\rmdefault}{\mddefault}{\updefault}
$\eta^2$
}}}
\put(2251,-1411){\makebox(0,0)[lb]{\smash{\SetFigFont{12}{14.4}{\rmdefault}{\mddefault}{\updefault}
1
}}}
\put(1126,-1561){\makebox(0,0)[lb]{\smash{\SetFigFont{12}{14.4}{\rmdefault}{\mddefault}{\updefault}
0
}}}
\put(3076,-1636){\makebox(0,0)[lb]{\smash{\SetFigFont{12}{14.4}{\rmdefault}{\mddefault}{\updefault}
-1
}}}
\put(4201,-1636){\makebox(0,0)[lb]{\smash{\SetFigFont{12}{14.4}{\rmdefault}{\mddefault}{\updefault}
0
}}}
\put(5326,-1636){\makebox(0,0)[lb]{\smash{\SetFigFont{12}{14.4}{\rmdefault}{\mddefault}{\updefault}
1
}}}
\put(676,-511){\makebox(0,0)[lb]{\smash{\SetFigFont{12}{14.4}{\rmdefault}{\mddefault}{\updefault}
$\eta$
}}}
\put(376,-2611){\makebox(0,0)[lb]{\smash{\SetFigFont{12}{14.4}{\rmdefault}{\mddefault}{\updefault}
}}}
\put(526,-2761){\makebox(0,0)[lb]{\smash{\SetFigFont{12}{14.4}{\rmdefault}{\mddefault}{\updefault}
\textbf{Figure 1}: $E = \{1, \eta, \eta^2 \}$
}}}
\put(3151,-2761){\makebox(0,0)[lb]{\smash{\SetFigFont{12}{14.4}{\rmdefault}{\mddefault}{\updefault}
\textbf{Figure 2}: $E = \{ -1,0,1 \}$
}}}
\end{picture}

\end{center}

\bigskip

There is an alternative way to describe the set $\tE$. To explain this, let $\Gamma$ be the Fuchsian group such that
$\Lambda (r) / \Gamma \cong \Omega$ 
and $g \circ \tau (z) = g(z)$ for all $\tau \in \Gamma$ and 
$z \in \Lambda (r)$. We denote by $\Gamma_\infty$ the orbit of $\infty$
under $\Gamma$ and by $d ( \cdot, \cdot )$ the hyperbolic distance between two points in $\Lambda (r)$. The set
\begin{equation}\label{diri}
 \D := \{ w \in \Lambda (r) : d (\infty, w) < d(\zeta, w) \; \mbox{for all } 
       \zeta \in \Gamma_\infty, \; \zeta \ne \infty \}
\end{equation}
is called the \emph{Dirichlet domain} associated with $\Gamma$ (or $g$) and centered at $\infty$. 
It is a fundamental region (cf. \cite{LE}, Section I-4).
Then it will be shown in Corollary~\ref{equi} that
$\tE = \partial g(\D)$, the boundary of the image of $\D$ under $g$.

Let $\U$ be the collection of compact sets in $\mathbb{C}$ 
with connected complements such that either $E$ is not connected, 
or $E$ is connected but $\partial E$ contains a singularity other than an outward cusp.

Our main result is:

\begin{theorem}\label{T}
\begin{enumerate}[(1)]
 \item $\tE$ is connected.
 \item If $E^\circ = \emptyset$, where $E^\circ$ denotes the interior of $E$,
       then the sequence $\{ \mu_k \}$ converges in the weak-star topology
       to the probability measure 
       \begin{equation}\label{mu}
       \mu (z) = \frac{1}{2 \pi} \Delta \log \delta (z)
       \end{equation}
       where 
       \begin{equation}\label{delta}
       \delta (z) = \limsup_{k \to \infty} |P_k (z)|^{1/k}, 
       \end{equation}
       and the support of $\mu$ is $\tE$. If, in addition, 
       $\overline{\mathbb{C}} \backslash E$ is simply-connected, then $\mu$ is
       the equilibrium distribution of the compact set $E$.
 \item Suppose $E^\circ$ is connected and $E \in \U$. 
       Then there is a subsequence of $\{ \mu_k \}$ that converges in
       the weak-star topology to the measure $\mu$ in \eqref{mu}, and its support is $\tE$.
       If, in addition, $\overline{\mathbb{C}} \backslash E$ is simply-connected, then $\mu$ is the
       equilibrium distribution of the compact set $E$.
\end{enumerate}
\end{theorem}

The polynomials $\{ P_k \}$ 
provide some useful tools in other branches of mathematics. For example,
A.~Atzmon, A.~Eremenko and M.~Sodin showed in \cite{AES} that if $a$ is an element
in a complex Banach algebra with unit and $E \subset \mathbb{C}$ is a compact set with
connected complement, then the spectrum of $a$ is included in $E$ if and only if
\[
 \limsup_{k \to \infty} \| P_k (a) \|^{1/k} \leq r,
\]
where $\{ P_k \}$ and $r = r(E)$ are as before. They also showed in the same paper that for a given analytic germ 
$f(z) = \sum_{k=0}^\infty f_k z^{-k-1}$ at infinity with values in a Banach space, 
the polynomials $\{ P_k \}$ can be used to determine
whether $f$ has analytic continuation to $\Omega = \overline{\mathbb{C}} \backslash E$ or not. 

For a formal Laurent series of the form \eqref{laurent},
it is known (cf. \cite{U3}) that the zeros of $P_k(z)$ are exactly the eigenvalues of the leading $k \times k$ 
principal submatrix of the infinite Toeplitz matrix
\[
  \begin{bmatrix}
     b_0    & b_1    & b_2    & b_3    & \cdots \\
      1     & b_0    & b_1    & b_2    & \cdots \\
      0     &  1     & b_0    & b_1    & \cdots \\
      0     &  0     &  1     & b_0    & \cdots \\
    \cdots  & \cdots & \cdots & \cdots & \cdots 
  \end{bmatrix}.
\]

\section{Known facts}\label{S2}
Suppose the Laurent expansion of $g$ is given by \eqref{laurent}, and we
assume that $\rho_0 := \limsup_{k \to \infty} |b_k|^{1/k} < \infty$
so that the series is convergent in $|w| > \rho_0$.
We also define $\delta_0$ as the smallest nonnegative number 
such that $g$ has a \emph{meromorphic} extension to 
$\Lambda(\delta_0) = \{ |w| > \delta_0 \} \cup \{ \infty \}$.
Now for every $z \in \mathbb{C}$, we write
$g^{-1} (z) := \{w \in \Lambda(\rho_0) : g(w) = z \}$ and 
$\tilde{g}^{-1} (z) := \{w \in \Lambda(\delta_0) : g(w) = z \}$.

Note that the multiplicity of a point 
$w \in \Lambda (\delta_0)$ is $m$ if $g'(w) = \cdots = g^{(m-1)} (w) = 0$ and $g^{(m)} (w) \ne 0$.
The following definition is due to J.~L.~Ullman (\cite{U1}, \cite{U2}).

\begin{definition}\label{C}
For every nonnegative integer $p$, $\tC_p$ (or $C_p$) is the set of all points 
$z \in \mathbb{C}$ such that
the points of largest absolute value in $\tgi (z)$ (or $g^{-1} (z)$, respectively) have total 
multiplicity $p$. 
\end{definition} 
From the definition, one can easily see that $\tC_0 = \mathbb{C} \backslash g(\Lambda(\delta_0))$ is compact and 
$\tC_1$ is an open set containing a neighborhood of infinity.

The next two statements are analogous to Theorem~1.3 and Theorem~1.4 of \cite{KS}, where
the theorems below are proved with $C_0, C_1, F_k$ and $\nu_k$ in place of $\tC_0, \tC_1, P_k$ and $\mu_k$, 
respectively.

\begin{theorem}[cf. \cite{KS}, Theorem~1.3]\label{T1}
If the interior of $\tC_0$ is empty, then the sequence $\{ \mu_k \}$ converges
in the weak-star topology to the measure 
$\mu$ in \eqref{mu} and
the support of $\mu$ is equal to $\partial \tC_1$.
If, in addition, $\overline{\mathbb{C}}=\tC_0 \cup \tC_1$, 
then $\mu$ is the equilibrium distribution of the compact set $\tC_0$. 
\end{theorem}

\begin{theorem}[cf. \cite{KS}, Theorem~1.4]\label{T2}
If the interior of $\tC_0$ is connected, then there is a subsequence of $\{ \mu_k \}$
that converges in the weak-star topology to the measure $\mu$ in \eqref{mu}, and the support
of $\mu$ is $\partial \tC_1$. If, in addition, 
$\overline{\mathbb{C}} = \tC_0 \cup \tC_1$, then $\mu$ is
the equilibrium distribution of the compact set $\tC_0$.
\end{theorem}

The proofs of Theorems~\ref{T1} and \ref{T2} are given in Appendix~\ref{appen}.
In fact, one may check that these theorems can be shown by
the same arguments for Theorems~1.3 and 1.4 in \cite{KS}, provided that Ullman's results
used in \cite{KS} are replaced by those listed below.

\begin{lemma}[\cite{U2}, Lemma~7; cf. \cite{KS}, Lemma~2.2]\label{LL}
Every $z_0 \in \tC_p$, $p \geq 2$, has a neighborhood $B(z_0, \epsilon) := \{ z : |z-z_0| < \epsilon \}$
such that $\partial \tC_1 \cap \overline{B(z_0, \epsilon)}$ consists of a finite number of analytic Jordan
arcs each joining $z_0$ to a point on the circle $\partial B(z_0, \epsilon)$. Any two arcs intersect only 
at $z_0$. The remaining points of $B(z_0, \epsilon)$ are in $\tC_1$.
\end{lemma}

From Lemma~\ref{LL}, we have an important corollary.

\begin{corollary}[cf. \cite{KS}, Corollary~2.3]\label{cor}
\[
\partial \tC_1 = \partial \tC_0 \cup \bigcup_{p \geq 2} \tC_p.
\]
\end{corollary}

A point $\xi \in \mathbb{C}$ is called 
a limit point of the zeros of $\{ P_k \}$ if there exist 
an increasing sequence $\{ k_j \}$ and a zero $\xi_j$ of $P_{k_j}$ for each $j$ such that 
$\xi = \lim_{j \to \infty} \xi_j$.

\begin{theorem}[\cite{U2}; cf. \cite{KS}, Theorem~1.2]\label{ta}
All limit points of the zeros of $\{ P_k \}$ are in 
$\mathbb{C} \backslash \tC_1$. Every boundary point of $\tC_1$ is a limit point of the zeros of $\{ P_k \}$.
\end{theorem}

\begin{proof}
The first statement is Lemma~3 of \cite{U2}, and the second statement is
given in the proof of Lemma~11 in the same paper.
\end{proof}

\begin{lemma}[\cite{U2}, Lemma~12]\label{connected}
The boundary of the unbounded component of $\tC_1$ is a connected set.
\end{lemma}

\begin{lemma}[\cite{U2}, Lemma~8; cf. \cite{KS}, Lemma~2.4]\label{24}
For every $\epsilon_0 > 0$ and $z_0 \in \tC_p$, $p \geq 2$, there are $\epsilon_1 > 0$ and
$z_1 \in \tC_q$, $q \geq 2$, such that
\begin{gather*}
B(z_1, \epsilon_1) \subset B(z_0, \epsilon_0), \\
B(z_1, \epsilon_1) \cap \tC_0 = \emptyset, \\
B(z_1, \epsilon_1) \cap \tC_1 = D_1 \cup D_2,
\end{gather*}
where $D_1$ and $D_2$ are disjoint non-empty domains. Moreover, there exist analytic 
functions $f_1$ and $f_2$ on $B(z_1, \epsilon_1)$ such that
\begin{align}
|f_1 (z)| > |f_2 (z)|, \quad & z \in D_1, \label{1}\\
|f_i (z)| = \delta (z), \quad & z \in D_i, \; i = 1,2, \label{2}
\end{align}
where $\delta (z)$ is defined in \eqref{delta}.
\end{lemma}

\begin{lemma}[\cite{U2}; cf. \cite{KS}, Lemma~3.1]\label{31} 
\begin{enumerate}[(a)]
\item For every $z \in \tC_0$,
 \begin{equation}\label{a}
 \delta (z) = \delta_0,
 \end{equation}
 where $\delta_0$ is defined in the first paragraph of this section.
\item For every $z \notin \tC_0$,
 \begin{equation}\label{b} 
 \delta (z) = \max \{ |w| : w \in \tgi (z) \}.
 \end{equation}
\item For every $z \in \tC_1$,
 \begin{equation}\label{c}
 \delta (z) = \lim_{k \to \infty} |P_k (z)|^{1/k};
 \end{equation}
 i.e., the $\limsup$ in \eqref{delta} can be replaced by $\lim$.
\end{enumerate}
\end{lemma}

An immediate corollary of \eqref{a} and \eqref{b} is:
\begin{corollary}[cf. \cite{KS}, Lemma~3.1]\label{cor2}
$\delta (z)$ is a continuous function on $\mathbb{C}$.
\end{corollary}

\section{Proof of Theorem~\ref{T}}
Suppose $E \subset \mathbb{C}$ is a compact set containing more than
two points such that its complement 
$\Omega = \overline{\mathbb{C}} \backslash E$ is connected, and let 
$g : \Lambda (r) \to \Omega$ be the uniformizing map with 
Laurent expansion \eqref{laurent} at infinity.

\begin{lemma}\label{corresp}
For given $z \in \Omega$, there is a one-to-one correspondence between
shortest curves (with respect to the hyperbolic metric in $\Omega$)
from $z$ to $\infty$ and points in $\tgi (z)$
with largest absolute value.
\end{lemma} 

\begin{proof}
Suppose $\gamma : [a, b] \to \Omega$ is a shortest curve from $z$ 
to infinity such that $\gamma(a) = z$ and $\gamma(b) = \infty$,
and let $\tilde{\gamma} : [a, b] \to \Lambda (r)$ be
the lifting curve of $\gamma$ 
such that $\tilde{\gamma}(b) = \infty$.
Because $g$ is a local isometry between $\Lambda(r)$ and $\Omega$ with hyperbolic
metrics on them, 
we know that $\tilde{\gamma}$ is a shortest curve in $\Lambda (r)$; i.e., 
the trace of $\tilde{\gamma}$ is the ray
\begin{equation}\label{ray}
\{ t w : 1 \leq t \leq \infty \}
\end{equation}
for $w = \tilde{\gamma} (a)$.
Now if $|w_0| >  |\tilde{\gamma} (a)|$ for some $w_0 \in \tgi (z)$,
let $\alpha$ be the ray \eqref{ray} with $w = w_0$. Then since 
$|w_0| >  |\tilde{\gamma} (a)|$, $\alpha$ is shorter than $\tilde{\gamma}$,
hence $g(\alpha)$ is shorter than $g(\tilde{\gamma}) = \gamma$.
Since $g(\alpha)$ connects $z$ to $\infty$, this contradicts to the choice of $\gamma$,
proving $|\tilde{\gamma} (a)| = \max \{ |w| : w \in \tgi (z) \}$. 
Conversely, if we have a point $w_0 \in \tgi (z)$ 
such that $|w_0| = \max \{ |w| : w \in \tgi (z) \}$,
then by the same argument above the image of the ray \eqref{ray} with $w=w_0$ is 
a shortest curve in $\Omega$ connecting $z=g(w_0)$ and $\infty$. 
Therefore the map $\gamma \mapsto \tilde{\gamma} (a)$ is bijective
between the set of shortest curves from $z$ to $\infty$ and the set
of points in $\tgi (z)$ with largest absolute value. The lemma follows.
\end{proof}

Note that $g'(w) \ne 0$ for any $w \in \Lambda (r)$ since $g$ is
a covering map. Therefore the following corollary is an immediate consequence of Lemma~\ref{corresp}.

\begin{corollary}\label{tE=C_p}
For each $z \in \Omega$, $z \in \tE$ if and only if $z \in \tC_p$ for some $p \geq 2$.
\end{corollary}

Let $\D$ be the Dirichlet domain associated with $g$
and centered at $\infty$, as we introduced in Section~\ref{intro}. 

\begin{corollary}\label{equi}
$\tE = \partial g(\D)$.
\end{corollary}

\begin{proof}
Because $\D$ is a fundamental region, $g(\overline{\D}) = \Omega$ and
$g |_\D$ is injective. 
Therefore $\partial E = \partial \Omega \subset \partial g(\D)$, 
$E^\circ \cap \partial g(\D) = \emptyset$, and
$\partial g(\D) \cap \Omega = g(\partial \D)$. Hence to prove the corollary, it suffices to
show that $\tE \cap \Omega = \partial g(\D) \cap \Omega = g (\partial \D)$.
 
Let $d(\cdot, \cdot)$ denote the hyperbolic distance between two points in $\Lambda (r)$,
and let $\Gamma$ be the corresponding Fuchsian group such that $\Lambda (r) / \Gamma \cong \Omega$.
Then for each $w_1, w_2 \in \Lambda(r)$ and $\tau \in \Gamma$, 
$d (w_1, w_2) = d(\tau(w_1), \tau (w_2))$. 
Therefore, $d(w, \infty) \leq d(w, \tau (\infty))$ for all $\tau \in \Gamma$
if any only if $d(w, \infty) \leq d( \tau(w), \infty)$ for all $\tau \in \Gamma$, or
$|w| \geq |\tau (w)|$ for all $\tau \in \Gamma$. This means that 
$w \in \overline{\D}$ if and only if
$|w| = \max \{ |\xi| : \xi \in \tgi ( g(w)) \}$. Now the corollary follows from Lemma~\ref{corresp}
because we have $w \in \partial \D$ if and only if there exist $w' \in \partial \D \backslash \{w \}$ 
and $\tau \in \Gamma$ such that $w = \tau (w')$, that is, 
$g(w) = g(w')$ (\cite{LE}, p.~37). 
\end{proof}

Note that the proof of this corollary also shows that $\Omega \backslash \tE = g(\D)$.

\begin{proof}[Proof of Theorem~\ref{T}]
Suppose $E \notin \U$. Then $\Omega$ is simply-connected, hence $\tE = \partial E = \partial \Omega$
is connected. If $E \in \U$, then as shown in \cite{KS} (p.~444), $g$ cannot be extended to $\Lambda(r_0)$ for
any $r_0 < r$; i.e., $r = \delta_0$, where $\delta_0$ is defined in the first paragraph
of Section~\ref{S2}. Therefore we have $E= \tC_0$, hence $\Omega \backslash \tE = \tC_1$
by Corollary~\ref{tE=C_p}; i.e., $\tC_1 = g( \D )$.
Now Lemma~\ref{connected} implies the statement (1) of Theorem~\ref{T},
since $\tC_1 = g( \D )$ is connected and $\tE = \partial g(\D)$.

To prove (2) of Theorem~\ref{T}, assume that  $E^\circ = \emptyset$.
If $E$ is not connected, then we have $\delta_0 = r$,
or $\tC_0 = E$, hence $\partial \tC_1 = \tE$ by Corollaries~\ref{cor} and \ref{tE=C_p}.
Since $\tC_0 = E$ has empty interior,
Theorem~\ref{T1} implies that $\{ \mu_k \}$ converges to $\mu$ and
its support is $\partial \tC_1 = \tE$. 
If $E$ is connected, we consider two cases: $\delta_0 < r$ and $\delta = r$. 
If $\delta_0 < r$, then every point $z \in E$ corresponds to at
least 2 points on $|w| = r$ counting muliplicity, hence $\tC_0 = \emptyset$, 
$\delta (z)$ is constant ($=r$) on $E$, and $E = \bigcup_{p \geq 2} \tC_p = \partial \tC_1$.
If $\delta_0 = r$, then $\delta (z)$ is constant ($=r=\delta_0$) on $E=\tC_0$ by Lemma~\ref{31}(a) and 
$E = \tE = \partial \tC_1$ as before. 
In any case $E = \tE = \partial \tC_1$ and $\delta (z)$ is constant on $E$.
Therefore by Theorem~\ref{T1}, 
the sequence $\{ \mu_k \}$ converges to $\mu = (2 \pi)^{-1} \Delta \log \delta$   
and its support is $E = \tE$. Furthermore, it is the equilibrium distribution of $E$ since $\delta$ is
constant on $E$. This proves (2) of Theorem~\ref{T}.

The statement (3) of Theorem~\ref{T} follows from Theorem~\ref{T2} by the same
argument as above. 
\end{proof}

One cannot replace $\mu_k$ in Theorem~\ref{T} by 
$\nu_k = (2 \pi k)^{-1} \Delta (\log |F_k|)$. To see this,
let $\eta = \exp (2 \pi i/3 )$ , $E = \{ 1, \eta, \eta^2 \}$, 
and $g : \Lambda (r) \to \Omega$ the universal covering map 
with Laurent expansion \eqref{laurent}. As in Example~\ref{example},
\[
\tE = \{ \eta^j t : j = 0,1,2 \mbox{ and } 0 \leq t \leq 1 \}.
\]  
Now we claim that there is a subsequence of $\{ \nu_k \}$ 
such that the support of its weak-star limit $\nu$ is
different from $\tE$ (Figure~4).
\bigskip

\begin{center}
 \begin{picture}(0,0)%
\includegraphics{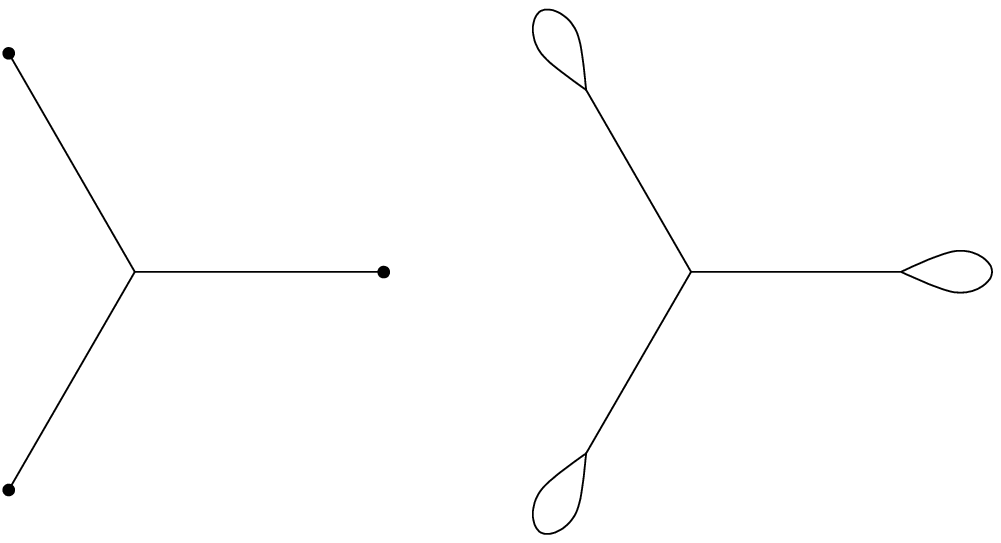}%
\end{picture}%
\setlength{\unitlength}{3947sp}%
\begingroup\makeatletter\ifx\SetFigFont\undefined%
\gdef\SetFigFont#1#2#3#4#5{%
  \reset@font\fontsize{#1}{#2pt}%
  \fontfamily{#3}\fontseries{#4}\fontshape{#5}%
  \selectfont}%
\fi\endgroup%
\begin{picture}(4774,2998)(601,-2669)
\put(751,-2086){\makebox(0,0)[lb]{\smash{\SetFigFont{12}{14.4}{\rmdefault}{\mddefault}{\updefault}
$\eta^2$
}}}
\put(2551,-1111){\makebox(0,0)[lb]{\smash{\SetFigFont{12}{14.4}{\rmdefault}{\mddefault}{\updefault}
1
}}}
\put(751, 14){\makebox(0,0)[lb]{\smash{\SetFigFont{12}{14.4}{\rmdefault}{\mddefault}{\updefault}
$\eta$
}}}
\put(1276,-1186){\makebox(0,0)[lb]{\smash{\SetFigFont{12}{14.4}{\rmdefault}{\mddefault}{\updefault}
0
}}}
\put(601,314){\makebox(0,0)[lb]{\smash{\SetFigFont{12}{14.4}{\rmdefault}{\mddefault}{\updefault}
}}}
\put(601,-2611){\makebox(0,0)[lb]{\smash{\SetFigFont{12}{14.4}{\rmdefault}{\mddefault}{\updefault}
\textbf{Figure 3}: The set $\tilde{E}$
}}}
\put(3151,-2611){\makebox(0,0)[lb]{\smash{\SetFigFont{12}{14.4}{\rmdefault}{\mddefault}{\updefault}
\textbf{Figure 4}: The support of $\nu$
}}}
\end{picture}

\end{center}
\bigskip

Let $\rho_0$ be the maximum of absolute values of the finite poles of $g$,
and let $\Omega_0 := g(\Lambda (\rho_0))$.
Then since $\rho_0 > r = r(E)$ and 
\begin{equation}\label{invrotation}
\eta^{-1} g (\eta w) = g(w) \; \mbox{for all } w \in \Lambda (r),
\end{equation} 
one can easily see that 
$\mathbb{C} \backslash \Omega_0 
= \mathcal{N}_0 \cup \mathcal{N}_1 \cup \mathcal{N}_2$,
where $\mathcal{N}_0$ is a closed neighborhood of 1 and 
$\mathcal{N}_j = \eta^j \mathcal{N}_0$, $j=1,2$.
Therefore from Definition~\ref{C} and
Corollary~2.3 of \cite{KS} (cf. Corollary~\ref{cor}),
\[
 \partial C_1 = \partial C_0 \cup \bigcup_{p \geq 2} C_p
              = \bigcup_{j=1}^3 \partial \mathcal{N}_j \cup \left(
                \tE \cap \Omega_0 \right),
\]
which is the set shown in Figure~4.

By Lemma~3.1 of \cite{KS} (cf. Lemma~\ref{31}), we have
$\limsup_{k \to \infty} |F_k (1)| = \rho_0$. Thus there exists a
subsequence $\{ F_{k_j} \}$ of $\{ F_k \}$ such that 
$\lim_{j \to \infty} |F_{k_j} (1)| = \rho_0$.
Considering \eqref{invrotation}, we also have
$\lim_{j \to \infty} |F_{k_j} (\eta)| = 
 \lim_{j \to \infty} |F_{k_j} (\eta^2)| = \rho_0$.
Now by Lemma~5.1 and the proof of Lemma~5.2 of \cite{KS}
(cf. Lemmas~\ref{51} and \ref{lll}), the sequence 
$\{ \nu_{k_j} \}$ converges in the weak-star topology to the measure 
\[
\nu (z) = \frac{1}{2 \pi} \Delta \log  \left( 
          \limsup_{k \to \infty} |F_k (z)|^{1/k} \right),
\]
and by Lemma~4.1 of \cite{KS} (cf. Lemma~\ref{4lemma}) 
the support of $\nu$ is $\partial C_1$. 
 
\section{Compact sets consisting of three points}\label{ex}
Let $\0 := \mathbb{C} \backslash \{0,1\}$ be a metric space equipped with
the hyperbolic metric on it, and we use the notations
$H^+ := \{ Im (z) > 0 \}$, $H^- := \{ Im (z) < 0 \}$,
$I_1 := \{ z \in \mathbb{R} : z < 0 \},$
$I_2 := \{ z \in \mathbb{R} : 0< z < 1 \},$  
$I_3 := \{ z \in \mathbb{R} : z >  1 \}$, and   
$I := I_1 \cup I_2 \cup I_3$.  
The following lemma is very trivial but plays a crucial role 
in our arguments later. 

\begin{lemma}\label{nopf}
Suppose $z_1 \in H^+$ and $z_2 \in H^+ \cup I$. Then there is a
unique shortest curve $\gamma$ connecting $z_1$ and $z_2$.
Moreover, the interior arc of $\gamma$ does not intersect $I$.
\end{lemma}

\begin{proof}
The proof is omitted here and left to the reader.
In fact, one may prove it by considering the symmetric property of $\0$.
\end{proof}

Let $\mathbb{D} := \{ z \in \mathbb{C} : |z| < 1 \}$ 
with the hyperbolic metric on it, and we assume that $G: \mathbb{D} \to \0$ is 
a holomorphic universal covering map such that
for a given point $a \in I_1 \cup H^+$, $G(0) = a$. 
Let  $\Gamma$ be the modular group on $\mathbb{D}$ such that 
$\mathbb{D} / \Gamma \cong \0$ 
and $G \circ \tau (z) = G(z)$ for all $\tau \in \Gamma$, 
and we denote by $\Gamma_0$ the orbit of the origin under $\Gamma$ 
and by $\D_0$ the Dirichlet domain with centered at the origin.
Note that each component of $G^{-1} (H^+)$ or $G^{-1} (H^-)$ is a
hyperbolic open triangle, and each side of such a triangle is
a geodesic curve and a component of $G^{-1} (I_j)$ for some $j=1,2,3$. 
Now we consider the cases $a \in I_1$ and $a \in H^+$ separately.

If $a \in I_1$, there are two hyperbolic triangles
$\triangle^+$ and $\triangle^-$ in $\mathbb{D}$
such that $G(\triangle^+) = H^+$, $G(\triangle^-) = H^-$ and 
$0 \in \overline{\triangle^+} \cap \overline{\triangle^-}$.
Suppose $w \in \triangle^+$. Then by Lemma~\ref{nopf},
there exists a unique shortest curve $\gamma : [t_0, t_1] \to H^+ \cup \{a\}$ 
such that $\gamma (t_0) = a$ and $\gamma (t_1) = G(w)$.
Let $\tilde{\gamma}$ be the lifting curve such that $\tilde{\gamma} (t_0) =0$.
Because $\gamma$ does not intersect $I$, 
$\tilde{\gamma}$ cannot intersect $G^{-1} (I)$, hence 
$\tilde{\gamma} (t_1) = w$. Now if $G(w') = G(w)$ for some $w' \ne w$,
the shortest curve $\alpha$ from $w'$ to $0$ intersects $G^{-1} (I)$,
thus $G(\alpha)$ intersects $I$. Therefore the length of $G(\alpha)$
is strictly greater than the length of $\gamma$, or
the the length of $\alpha$ is strictly greater than the length of $\tilde{\gamma}$.
This shows that $w \in \D_0$. Because a similar argument holds for any $w \in \triangle^-$,
we have
\begin{equation}\label{DD}
 \triangle^+ \cup \triangle^- \subset \D_0.
\end{equation}
Since $\D_0$ is a fundamental region, $G$ is univalent in $\D_0$ and 
$G(\overline{\D}_0) = \0$.
Thus \eqref{DD} in fact shows that 
\[
\D_0 = \left( \overline{\triangle^+} \cup \overline{\triangle^-} \right)^\circ.
\]
Note that in this case $G(\partial \D_0) = I_2 \cup I_3$.

We next consider the case $a \in H^+$.
Let $\triangle_0$ be the triangle such that $0 \in \triangle_0$ and 
$G(\triangle_0) = H^+$.
For each $j \in \{ 1,2,3 \}$, we denote by $L_j$ 
the side of $\triangle_0$ such that $G(L_j) = I_j$,
and let $\triangle_j$ be the hyperbolic triangle which is obtained by
reflecting $\triangle_0$ with respect to $L_j$. Similarly,
we denote by $\triangle_{j,k}$ the triangle obtained by reflecting
$\triangle_j$ through the side over $I_k$, $k \ne j$. Finally,
let $\zeta_{j,k}$ be the point in $\triangle_{j,k}$ such that 
$G(\zeta_{j,k}) = a$. See Figure~5.

\bigskip

\begin{center}
 \begin{picture}(0,0)%
\includegraphics{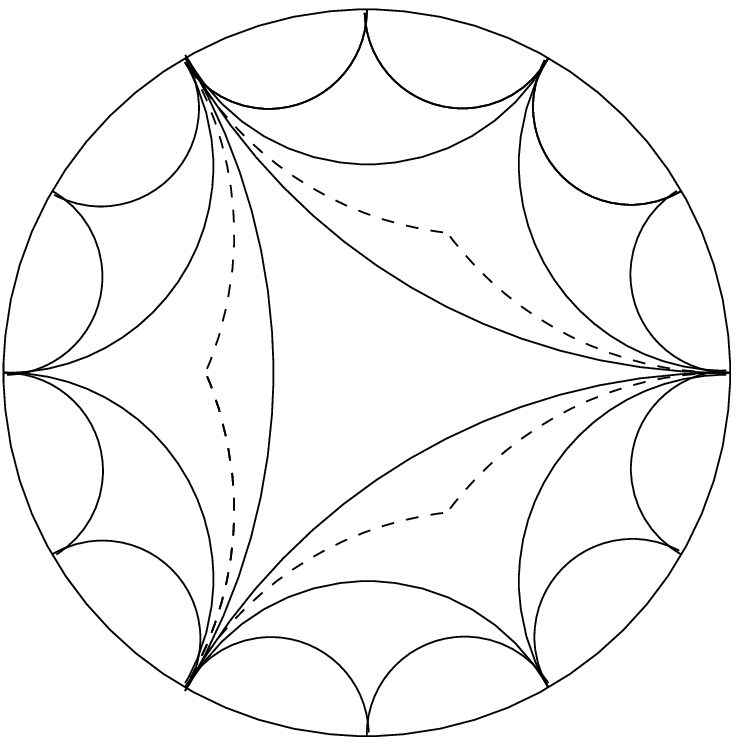}%
\end{picture}%
\setlength{\unitlength}{3947sp}%
\begingroup\makeatletter\ifx\SetFigFont\undefined%
\gdef\SetFigFont#1#2#3#4#5{%
  \reset@font\fontsize{#1}{#2pt}%
  \fontfamily{#3}\fontseries{#4}\fontshape{#5}%
  \selectfont}%
\fi\endgroup%
\begin{picture}(3504,4005)(1343,-4094)
\put(3151,-2311){\makebox(0,0)[lb]{\smash{\SetFigFont{12}{14.4}{\rmdefault}{\mddefault}{\updefault}
$L_1$
}}}
\put(3751,-2236){\makebox(0,0)[lb]{\smash{\SetFigFont{12}{14.4}{\rmdefault}{\mddefault}{\updefault}
$\ell_{1,2}$
}}}
\put(2926,-2686){\makebox(0,0)[lb]{\smash{\SetFigFont{12}{14.4}{\rmdefault}{\mddefault}{\updefault}
$\ell_{1,3}$
}}}
\put(3451,-2761){\makebox(0,0)[lb]{\smash{\SetFigFont{12}{14.4}{\rmdefault}{\mddefault}{\updefault}
$\triangle_1$
}}}
\put(2851,-3061){\makebox(0,0)[lb]{\smash{\SetFigFont{12}{14.4}{\rmdefault}{\mddefault}{\updefault}
$\triangle_{1,3}$
}}}
\put(3976,-2536){\makebox(0,0)[lb]{\smash{\SetFigFont{12}{14.4}{\rmdefault}{\mddefault}{\updefault}
$\triangle_{1,2}$
}}}
\put(3076,-1411){\makebox(0,0)[lb]{\smash{\SetFigFont{12}{14.4}{\rmdefault}{\mddefault}{\updefault}
$L_2$
}}}
\put(2551,-2161){\makebox(0,0)[lb]{\smash{\SetFigFont{12}{14.4}{\rmdefault}{\mddefault}{\updefault}
$L_3$
}}}
\put(4126,-1261){\makebox(0,0)[lb]{\smash{\SetFigFont{12}{14.4}{\rmdefault}{\mddefault}{\updefault}
$\cdot \; \zeta_{2,1}$
}}}
\put(3076,-1936){\makebox(0,0)[lb]{\smash{\SetFigFont{12}{14.4}{\rmdefault}{\mddefault}{\updefault}
$\cdot \; 0$
}}}
\put(2851,-1711){\makebox(0,0)[lb]{\smash{\SetFigFont{12}{14.4}{\rmdefault}{\mddefault}{\updefault}
$\triangle_0$
}}}
\put(3076,-661){\makebox(0,0)[lb]{\smash{\SetFigFont{12}{14.4}{\rmdefault}{\mddefault}{\updefault}
$\cdot \; \zeta_{2,3}$
}}}
\put(2026,-2536){\makebox(0,0)[lb]{\smash{\SetFigFont{12}{14.4}{\rmdefault}{\mddefault}{\updefault}
$\cdot \; \zeta_{3,1}$
}}}
\put(2026,-4036){\makebox(0,0)[lb]{\smash{\SetFigFont{12}{14.4}{\rmdefault}{\mddefault}{\updefault}
\textbf{Figure 5}: The case $a \in H^+$
}}}
\end{picture}

\end{center}

\bigskip

By Lemma~\ref{nopf} and the same argument as above, 
$\overline{\triangle}_0 \subset \D_0$. Similarly, the closed triangle
$\overline{\triangle}_{j,k}$ is contained in the Dirichlet domain
with center at $\zeta_{j,k}$. Since $\D_0$ is connected, we have
\begin{equation}\label{zeta}
 \overline{\triangle}_0 \subset \D_0 \subset \bigcup_{j=0}^4 \triangle_j \cup \bigcup_{j=1}^3 L_j.
\end{equation}

Let $A := \{ \zeta_{j,k} : 1 \leq j, k \leq 3, \; j \ne k \}$.
Then for all $w \in \D_0$ and 
$\zeta \in \Gamma_0 \backslash (\{ 0 \} \cup A)$, 
the shortest curve connecting these two points 
must pass through a point $w' \in \triangle_{j,k}$ for some $j$ and $k$, 
$j \ne k$. Therefore, denoting by $d_\mathbb{D}$ the hyperbolic distance in $\mathbb{D}$, we have 
\[
d_\mathbb{D} (w, \zeta) = d_\mathbb{D} (w, w') + d_\mathbb{D} (w', \zeta) > d_\mathbb{D} (w, w') + d_\mathbb{D} (w', \zeta_{j,k})
\geq d_\mathbb{D} (w, \zeta_{j,k}),
\]
because $w' \in \triangle_{j,k}$ is contained in the Dirichlet domain 
with center at $\zeta_{j,k}$. This implies that
\begin{equation}\label{6d}
\D_0 = \bigcap_{\zeta \in A} 
      \{ w : d_\mathbb{D} (w, 0) <  d_\mathbb{D} (w, \zeta) \}.
\end{equation}

We obtained $\triangle_{j,k}$ by reflecting $\triangle_0$
through a side over $I_j$ and then a side over $I_k$. 
The triangle $\triangle_0$ can be obtained from $\triangle_{k,j}$ by the same way.
Therefore, there is $\tau \in \Gamma$ such that $\tau (0) = \zeta_{j,k}$
and $\tau (\zeta_{k,j}) = 0$. Then with the notation
\[
\ell_{j,k} := \{ w : d_\mathbb{D} (w, 0) =  d_\mathbb{D} (w, \zeta_{j,k}) \},
\]
we have $\tau (\ell_{k,j}) = \ell_{j,k}$; i.e., $G(\ell_{j,k}) = G(\ell_{k,j})$.
Also note that $\ell_{j,k}$ is a geodesic curve which separates 
the two sides of $\triangle_j$ lying over $I_j$ and $I_k$,
because $\overline{\triangle}_0$ and $\overline{\triangle}_{j,k}$ are
contained in the Dirichlet domains with centers at $0$ and $\zeta_{j,k}$, respectively.
Since each sides of $\triangle_j$ are also geodesic, we conclude that 
one end of $\ell_{j,k}$ approaches to the common vertex
(at infinity) of $\triangle_0$ and $\triangle_{j,k}$, and $\ell_{j,k}$
intersects the side of $\triangle_j$ which is over $I_l$, $l \ne j, k$.
In particular, this implies that $\ell_{j,k} \cap \ell_{j,k'} \ne \emptyset$ for $k \ne k'$.

Now one can easily see that $\D_0$ is a hexagon with exactly 3 vertices at 
infinity, which are in fact the vertices of $\triangle_0$. 
Moreover along $\partial \D_0$, finite and infinite vertices are placed
alternatively and the $G$-images of the two sides sharing  
a common infinite vertex are same. Therefore, the three finite vertices
are mapped to the same point, say $b$, and $G(\partial \D_0)$ is a 
tripod with center at $b$ such that each leg of it 
is a hyperbolic geodesic curve connecting $b$ to
one of the points $0,1,\infty$. 
Now we are ready to prove the following proposition.

\begin{proposition}\label{prop}
Suppose $E$ consists of three points. If the points in $E$ are on a 
straight line, $\tE$ is the line segment connecting points in $E$; i.e.,
$\tE$ is the convex hull of $E$. Otherwise $\tE$ is a tripod.
\end{proposition}

\begin{proof}
Let $g: \Lambda (r) \to \overline{\mathbb{C}} \backslash E$ 
be a holomorphic universal map. We choose a linear
transformation $T$ such that $T(E) = \{0,1, \infty \}$ and 
$a:=T(\infty) \in H^+ \cup I_1$. Since every linear transformation
maps circles onto circles, we see that $a \in I_1$ if and only if
$\infty$ is on the circle passing through the points in $E$; i.e.,
if and only if the points in $E$ are on a straight line.
Now let $G(w) := T \circ g (r/w)$. Then $G$ is a holomorphic 
universal covering map from $\mathbb{D}$ to $\0$ such that $G(0) = a$.
Now the proposition follows from Corollary~\ref{equi} and the arguments
preceding the proposition, since $T$ is a conformal map
and the map $w \mapsto r/w$ sends 
the Dirichlet domain in $\mathbb{D}$ with center at the origin
onto the Dirichlet domain in $\Lambda(r)$ with center at infinity.
\end{proof}
 
\appendix
\section{Appendix: Proofs of Theorems~\ref{T1} and \ref{T2}}\label{appen}
For a measure $\omega$ on $\mathbb{C}$, we denote its logarithmic potential by
\[
 \PP_\omega (z) := - \int_{\mathbb{C}} \log |z-\xi| d\omega(\xi)
\]
and let $\PP(z) := \log \delta (z) = \limsup_{k \to \infty} k^{-1} \log |P_k (z)|$. 
Note that $\mu (z) = (2 \pi)^{-1} \Delta \PP(z)$ by the definition of $\PP$ and \eqref{mu}.

\begin{lemma}\label{4lemma}
\begin{enumerate}[(a)]
 \item $\PP$ is subharmonic on $\mathbb{C}$.
 \item $\PP$ is harmonic on $\tC_1 \cup (\tC_0)^\circ$, but not at points of 
       $\partial \tC_1$.
 \item $\mu$ is a probability measure with support $\partial \tC_1$. 
 \item $\PP_\mu (z) = - \PP (z)$ for all $z \in \mathbb{C}$. 
\end{enumerate}
\end{lemma}

\begin{proof}
Note that the limit superior of a sequence of subharmonic functions is subharmonic
if it is upper semi-continuous. Since $\PP$ is upper semi-continuous by
Corollary~\ref{cor2} and $k^{-1} \log |P_k|$ is subharmonic for all $k$, (a) follows.

By \eqref{a}, $\PP$ is constant on $\tC_0$ hence harmonic on $(\tC_0)^\circ$.
If $z \in \tC_1$, there exist a neighborhood $N$ of $z$
and an inverse branch of $g$, say $f$, defined on $N$ such that
\[
 |f(\xi)| = \max \{ |\zeta|: \zeta \in \tgi (\xi) \}, \quad \mbox{for } \xi \in N.
\]
Therefore by \eqref{b}, $\PP(\xi) = \log |f(\xi)|$ in $N$ 
hence harmonic since $|f (\xi)| > \delta_0 \geq 0$.

The function $\PP$ is not harmonic on $\partial \tC_0$ 
since Lemma~\ref{31} implies that $\PP (z)=\delta_0$ 
for all $z \in \tC_0$ but $\PP (z) > \delta_0$ if $z \notin \tC_0$. 
We next show that $\PP$ is not harmonic at a point $z \in \tC_p$, $p \geq 2$.
If it is not the case, there exists a neighborhood $N$ of $z$ such that
$\PP$ is harmonic on $N$. 
Then by Lemma~\ref{24}, there exist a subdomain $N' \subset N$,
two disjoint domains $D_1, D_2$ such that $D_1 \cup D_2 = N' \cap \tC_1$,
and analytic functions $f_1$ and $f_2$ satisfying 
\eqref{1} and \eqref{2}. This implies $\PP(\xi) - \log |f_2 (\xi)|$ is a harmonic
function which is positive on $D_1$ and zero on $D_2$, which is 
impossible. Therefore $\PP$ is not harmonic on $\bigcup_{p \geq 2} \tC_p$.
Now (b) follows from Corollary~\ref{cor}.

Since $\PP$ is subharmonic by (a), $\mu = (2 \pi)^{-1} \Delta \PP$ is
a measure. Moreover by (b), the support of $\mu$ is $\partial \tC_1$.
Therefore to show (c), it suffices to show that $\mu (\mathbb{C}) =1$.
Note that since $\mu$ has compact support and $\PP$ is harmonic off 
the set $\partial \tC_1$,
the Riesz Decomposition Theorem (\cite{LA}, Theorem~II.21) implies that  
\begin{equation}\label{u}
u(z) := \PP_\mu (z) + \PP (z)
\end{equation}
is harmonic on any bounded domain
$D$ containing $\partial \tC_1$. Because this is true for any arbitrary large
domain $D$, $u$ is in fact harmonic on $\mathbb{C}$, hence constant. 

Let $f$ be the inverse of $g$ defined on a neighborhood of $\infty$
such that $f(\infty)=\infty$. Then it is easy to see from \eqref{laurent} that
$f(z)  = z + O(1)$ as $z \to \infty$. Moreover by \eqref{b}, 
$\delta (z) = |f(z)|$ for sufficiently large $z$.
Thus $\PP(z) = \log |f(z)| = \log |z| + o(1)$ as $z \to \infty$. 
But since $\PP_\omega (z) = - \omega (\mathbb{C}) \log |z| + o(1)$ for any 
finite measure $\omega$ with compact support, we conclude from \eqref{u} that 
$\mu (\mathbb{C}) = 1$ and $u \equiv 0$, which shows (c) and (d) 
simultaneously. This completes the proof.
\end{proof}

Recall that $\mu_k (z) = (2\pi k)^{-1} \Delta (\log |P_k (z)|)$. Since $P_k$ is a monic
polynomial of degree $k$, it can be shown that $\PP_{\mu_k} = - k^{-1} \log |P_k|$
by the same argument for Lemma~\ref{4lemma}(d).

\begin{lemma}\label{51}
Let $\omega$ be any weak-star limit of $\{\mu_k\}$. Then
\begin{align}
 \PP_{\omega} (z) = \PP_\mu (z), & \quad z \in \tC_1 \cup \partial \tC_1, \label{P=P}\\
 \PP_{\omega} (z) \geq \PP_\mu (z), & \quad z \in \mathbb{C}. \label{P>P}
\end{align} 
\end{lemma}
 
\begin{proof}
Suppose $\{ \mu_{k_j} \}$ converges to $\omega$ in the weak-star topology.
By \eqref{c} and Lemma~\ref{4lemma}(d), 
\[
  \lim_{j \to \infty} \PP_{\mu_{k_j}} (z) 
  = - \lim_{j \to \infty} k_j^{-1} \log |P_{k_j} (z)| 
  = - \log \delta (z) = - \PP (z) = \PP_\mu (z)
\]
for all $z \in \tC_1$.
Therefore the Lower Envelope Theorem (\cite{LA}, Theorem~3.8) implies that
$\PP_\mu (z) = \PP_\omega (z)$ for all $z \in \tC_1$ except on a set
of logarithmic capacity zero. On the other hand, Theorem~\ref{ta}
implies that the support of  
$\PP_\omega$ is contained in $\mathbb{C} \backslash \tC_1$,
hence it is harmonic in $\tC_1$. Because $\PP_\mu = - \PP(z)$ is also
harmonic in $\tC_1$ (Lemma~\ref{4lemma}(b)), we conclude that  
$\PP_\mu (z) =  \PP_\omega (z)$ for all $z \in \tC_1$.

Corollary~\ref{cor2} implies that $\PP_\mu (z) = - \log \delta (z)$ is 
continuous if $\delta(z) >0$. If $\delta (z) = 0$ (this happens only when $z \in \tC_0$
and $\delta_0 = 0$), $\PP_\mu (z) = \infty$.
Therefore we have
\begin{equation}\label{eq}
\PP_\omega (z) \leq \PP_\mu (z), \quad \mbox{for all } z \in \partial \tC_1,
\end{equation}
because $\PP_\omega$ is lower semi-continuous and $\PP_\omega = \PP_\mu$ on $\tC_1$.

If $z \in \tC_1$ and $z$ approaches to $\partial \tC_0$, we have
$\PP_\omega (z) = \PP_\mu (z) \to - \log \delta_0$. Therefore,
the minimum principle implies that 
$\PP_\omega (z) \geq - \log \delta_0 = \PP_\mu (z)$ for all $z \in \tC_0$. 
Note that combining this result with \eqref{eq}, we also have 
$\PP_{\omega} (z) = \PP_\mu (z)$ for all $z \in \partial \tC_0$.

Now it remains to show that $\PP_{\omega} (z) \geq \PP_\mu (z)$ for all 
$z \in \tC_p$, $p \geq 2$. 
But for sufficiently small $\epsilon$, Lemma~\ref{LL} implies that 
the circle $\{ \xi : |z - \xi| = \epsilon \}$ is contained in $\tC_1$ 
except finitely many points. Since $\PP_\omega$ is
superharmonic and $\PP_\omega (\xi) = \PP_\mu (\xi)$ for  $\xi \in \tC_1$,
\[
 \PP_\omega (z) \geq 
    \frac{1}{2 \pi \epsilon} \int_{|z - \xi| = \epsilon} \PP_\omega (\xi) |d\xi|
    = \frac{1}{2 \pi \epsilon} \int_{|z - \xi| = \epsilon} \PP_\mu (\xi) |d\xi|.
\]
By letting $\epsilon \to 0$, we get $\PP_\omega (z) \geq \PP_\mu (z)$ since
$\PP_\mu (\xi)$ is continuous at $z \in \tC_p$, $p \geq 2$. This completes the proof.
\end{proof}

\begin{lemma}\label{lll}
Suppose $(\tC_0)^\circ \ne \emptyset$ and $U$ is a component of $(\tC_0)^\circ$.
Then there exists a subsequence of $\{ \mu_k \}$ that converges in the 
weak-star topology to a measure $\omega$ such that
\[
 \PP_\omega (z) = \PP_\mu (z) \quad \mbox{for all } z \in \overline{U}.
\]
\end{lemma}

\begin{proof}
Pick a point $z_0 \in U$. By \eqref{delta} there is a subsequence $\{ \mu_{k_j} \}$
such that
\[
 \lim_{j \to \infty} \PP_{\mu_{k_j}} (z_0)
 = - \lim_{j \to \infty} k_j^{-1} \log |P_{k_j} (z_0)| = - \log \delta_0.
\]
Now let $\omega$ be a weak-star limit of a subsequence of $\{ \mu_{k_j} \}$. 
Then by the Principle of Descent (\cite{LA}, Theorem~1.3), we have
$\PP_\omega (z_0) \leq - \log \delta_0$. Since \eqref{P=P} says that 
$\PP_\omega (z) =-\log \delta_0$ for all $z \in \partial U$, the minimum principle implies that 
$\PP_\omega = - \log \delta_0 = \PP_\mu$ in $U$. 
\end{proof}

\begin{proof}[Proofs of Theorems~\ref{T1} and \ref{T2}]
Suppose $(\tC_0)^\circ$ is empty and let $\omega$ be a weak-star limit of 
any convergent subsequence of $\{ \mu_k \}$. Then by \eqref{P=P}, we have $\PP_\omega (z) = \PP_\mu (z)$
for all $z \in \mathbb{C}$. 
Therefore $\omega = - (2\pi)^{-1} \Delta \PP_\omega  
=  - (2\pi)^{-1} \Delta \PP_\mu  = \mu$ and the first statement of 
Theorem~\ref{T1} follows from Lemma~\ref{4lemma}. If, in addition, 
$\mathbb{C} = \tC_0 \cup \tC_1$, then $\mu$ is the equilibrium distribution
of $\tC_0$ since the support of $\mu$ is $\partial \tC_1 = \partial \tC_0$ and
$\PP_\mu$ is constant on $\tC_0$ (\cite{TS}, Theorem~III.15).
This completes the proof of Theorem~\ref{T1}.

Now suppose $(\tC_0)^\circ$ is connected. Then by Lemma~\ref{lll}, 
there exists a subsequence $\{ \mu_{k_j} \}$ of $\{ \mu_k \}$ converging
to a measure $\omega$ such that $\PP_\omega = \PP_\mu$ in $\tC_0$.
Since $\PP_\omega (z) = \PP_\mu (z)$ for all $z \in \mathbb{C} \backslash \tC_0$ 
by \eqref{P=P}, this shows that $\PP_\omega = \PP_\mu$ in $\mathbb{C}$
hence $\omega = \mu$. Now Theorem~\ref{T2} follows from Lemma~\ref{4lemma}.
\end{proof}

\end{document}